\documentclass[centertags,12pt]{amsart}
\usepackage{latexsym}
\usepackage{amssymb}

\numberwithin{equation}{section}
\makeatletter
\renewcommand{\subsection}{\@stratsection
{subsection}{2}{0mm}{\baselineskip}{-0.25cm}
{\normalfont\normalize\bf}}
\makeatother

\newtheorem{prop}{Proposition}[section]
\newtheorem{thm}[prop]{Theorem}
\newtheorem{lemma}[prop]{Lemma}
\newtheorem{cor}[prop]{Corollary}
\theoremstyle{definition}
\newtheorem{defi}[prop]{Definition}
\newtheorem{rem}[prop]{Remark}
\newtheorem{ex}[prop]{Example}

{\theoremstyle{remark}
\newtheorem*{claim*}{Claim}}

\newcommand{\scrp }{{\mathcal P}}

\newcommand{\scrt }{{\mathcal T}}
\newcommand{\scrx }{{\mathcal X}}

\newcommand{\numbers }{{\mathbb N}}

\newcommand{\no }{{\mathbb N}_0}

\newcommand{\ff }{{\mathbb F}}

\newcommand{\fq }{{\mathbb F}_q}

\newcommand{\bc }{{\bf c}}

\newcommand{\bh }{{\bf h}}

\newcommand{\Hline}{\noindent\underline{\hspace{15.8cm}}}

\begin{document}

\author[A.~Campillo]{Antonio Campillo}
\author[J.I.~Farran]{Jose Ignacio Farran}
\author[C.~Munuera]{Carlos Munuera}

\title[AG-codes related to Arf semigroups]{On the parameters of algebraic geometry \\ codes related to Arf semigroups}
\address{A. Campillo, Departamento de Algebra y Geometria, Fac. de Ciencias, Universidad de Valladolid, 
Prado de la Magdalena SN, 47005 Valladolid, Castilla, Spain}
\email{campillo@agt.uva.es}
\address{J.I. Farran, Departamento de Matematica Aplicada, ETSII, Universidad de Valladolid, 
Paseo del Cauce SN, 47011 Valladolid, Castilla, Spain}
\email{ignfar@eis.uva.es}
\address{C. Munuera, Departamento de Matematica Aplicada, ETSA, Universidad de Valladolid, 
Avda. Salamanca SN, 47014 Valladolid, Castilla, Spain}
\email{cmunuera@modulor.arq.uva.es}

\thanks{Date: October 25, 1999}

\begin{abstract} 
In this paper we compute the order (or Feng-Rao) bound on the minimum distance of one-point algebraic geometry codes 
$C_{\Omega}(\scrp,\rho_lQ)$, when the Weierstrass semigroup at the point $Q$ is an Arf semigroup. 
The results developed to that purpose also provide the dimension of the improved geometric Goppa codes 
related to these $C_{\Omega}(\scrp,\rho_lQ)$.\\

\end{abstract}
\maketitle

\section{Introduction}

Let $\fq$ be a finite field and $F$ a function field over $\fq$. The construction of algebraic geometry 
(or geometric Goppa) codes from $F$ is well known (see \cite{HPvL}). Take a rational place $Q$ and let 
$K_{\infty}(Q)$ be the set (ring) of functions having no poles outside $Q$. Let $S=S(Q)$ be the Weierstrass 
semigroup  of $Q$, that is $S=\{ -v_Q(f) \ | \  f\in K_{\infty}(Q)\}$, where $v_Q$ is the valuation at $Q$. 
Usually we shall write $S$ as an enumeration of its elements in increasing order, $S=\{\rho_1=0<\rho_2<\cdots\}$. 
For a positive integer $m$, we also consider $L(mQ)=\{ f\in K_{\infty}(Q) \ | \  v_Q(f)\ge -m\}$. Given a set of 
$n$ distinct rational places in $F$,  $\scrp=\{ P_1,\cdots,P_n\}$, such that $Q\not\in \scrp$, we consider 
the evaluation map
$$
ev_{\scrp}:K_{\infty}(Q)\longrightarrow \fq^n \; , \;
ev_{\scrp}(f)=(f(P_1),\cdots,f(P_n))
$$
and define the (one-point) algebraic geometry code $C_{\Omega}(\scrp,\rho_lQ)=ev_{\scrp}(L(\rho_lQ))^{\perp}$, 
that is, if for $i=1,2,\cdots$, we choose a function $h_i\in K_{\infty}(Q)$ such that $-v_Q(h_i)=\rho_i$, 
then $C_{\Omega}(\scrp,\rho_lQ)$ is defined by the system of parity checks $\bh_1,\cdots,\bh_l$, 
with $\bh_j=ev_{\scrp}(h_j)$. For simplicity, from now on, we shall write $C_l$ instead of $C_{\Omega}(\scrp,\rho_lQ)$ 
if no confusion arises.

The parameters of $C_l$ are as follows: its length is obviously $n$ and its dimension is at least $n-l$, 
with equality if $\rho_l<n$. When $\rho_l\ge n$, then some of the checks $\bh_1,\cdots,\bh_l$ can be dependent, 
and the exact value of the dimension can be computed with the help of the Riemann-Roch theorem. Thus, 
these two parameters are, at least theoretically, easy to compute. On the contrary, 
the minimum distance of $C_l$ is often difficult to compute. A general lower bound on $d(C_l)$ 
is given by the Goppa bound (or {\em Goppa designed minimum distance}), $d(C_l)\ge d_{G}(l)=l+1-g$, 
where $g$ is the genus of $F$ (that is, the number of gaps in $S$). A better bound is the so-called {\em Feng-Rao} 
or {\em order} bound $d_{ORD}(C_l)$, defined as follows (see \cite{FRB} and \cite{HPvL}): 
for a pole $\rho\in S$, let us consider the set
$$
A[\rho]=\{ p\in S | \rho-p\in S\},
$$
Then, the {\em order bound} on the minimum distance of $C_l$ is
$$
d_{ORD}(l)=\min\{ \# A[\rho] | \rho\in S, \rho\ge\rho_{l+1} \}
$$
and it holds that $d_G(l)\le d_{ORD}(l)\le d(C_l)$. A remarkable property of the order bound is that it is computed 
only in terms of the semigroup $S$ (that is, without any relation neither with $F$ nor the set $\scrp$).

A way to improve algebraic geometry codes was introduced by Feng and Rao in \cite{FR}.  
For a positive integer $d$ let us consider the set
$$
R_d=\{ i \ | \ \# A[\rho_i]<d\} .
$$
The {\em improved geometric Goppa code} $\tilde{C}(d)$ is defined as
$$
\tilde{C}(d)=\{ \bc\in\ff^n \ | \ \bc\cdot\bh_i=0 \mbox{ for all } i\in R_d\}
$$
(see \cite{FR,HPvL,PT}). The parameters of  $\tilde{C}(d)$ are as follows: its minimum distance is at least $d$. 
Furthermore, if $d=d_{ORD}(l)$, then $C_l\subseteq\tilde{C} (d)$ (which explains the meaning of the term \lq improved'), 
so  its dimension is at least the dimension of $C_l$. On the other hand, it is clear that  
$\dim \tilde{C}(d)\ge n-\# R_d$, with equality if $2c\le n$ and $1\le d\le 2r-2$, where $c=\rho_r$ 
is the conductor of $S$. Thus, here we find again that the unknown parameters of the code can be estimated 
in terms of the semigroup $S$ (and more precisely, they are closely related to the A-sets $A[\rho]$ in $S$).

In this paper, we show how to compute both, the order bound on the minimum distance of an (one-point) 
algebraic geometry code and the redundancy of the corresponding improved code, when the involved semigroup 
$S$ is an Arf semigroup. The organization of the paper is as follows: Arf semigroups, their main properties 
and some examples are presented in section 2. In section 3, we show how to deal with the sets $A[\rho]$ for 
Arf semigroups. These results are used in section 4 for computing the order bound on the minimum distance 
of $C_l$, and again in section 5 in order to give a formula for $\# R_d$. With regard to this last section, 
we have to point out that  the study of the sequence $(\# R_d)$ has been already treated in the paper \cite{PT} 
by Pellikaan and Torres, and in fact, our section 5 can be viewed as a continuation of that paper. In particular, 
we simplify and extend some results stated there, and solve some open problems from it.

\section{Arf semigroups}

Let $S=\{\rho_1=0<\rho_2<\cdots\}$ be a numerical semigroup. Let $c=\rho_r$ be the conductor of $S$ 
and let $g=c-r+1$ be its genus. The elements $\rho\in S$ will be called {\em poles} and the elements 
$n\in\no\setminus S$ will be called {\em gaps}.

\begin{defi} \label{definicion}
$S$ is called an {\em Arf semigroup} if for every $i,j,k\in\numbers$ with 
$i\ge j\ge k$, it holds that $\rho_i+\rho_j-\rho_k\in S$. 
\end{defi}

Arf semigroups were introduced by C. Arf in \cite{Ar} as the semigroups of values of the so called 
{\em Arf one-dimensional local rings}, which geometrically correspond to curve singularities being maximal 
among the class of singularities with the same resolution type (see \cite{Li} for details).

\begin{rem} If $\rho_i\ge c$, then for every $j,k$, with $i\ge j\ge k$, we have $\rho_i+\rho_j-\rho_k\in S$. 
Thus, the condition stated in definition \ref{definicion} should be imposed only in the range that $k\le j\le i<r$.\par
\end{rem}

The defining Arf condition can be changed by another, apparently weaker, property.

\begin{prop} \label{debil}
Let $S$ be a semigroup. The following conditions are equivalent:
\newline a) $S$ is Arf;
\newline b) for every two positive integers $i,k$, with $i\ge k$, it holds that $2\rho_i-\rho_k\in S$.
\end{prop}
\begin{proof}
Obviously every Arf semigroup verifies b). Conversely, let us assume b) and let $i,j,k$, be positive integers 
such that $k\le j\le i<r$. We have to prove that $m=\rho_i+\rho_j-\rho_k\in S$. If $i=j$ or $j=k$, then it is clear. 
Otherwise, if $k<j<i$, let $i_0=i,j_0=j,k_0=k$, and write
$$
m=\rho_{i_0}+\rho_{j_0}-\rho_{k_0}=(2\rho_{j_0}-\rho_{k_0})+\rho_{i_0}-\rho_{j_0}.
$$
Note that $2\rho_{j_0}-\rho_{k_0}\in S$ and $2\rho_{j_0}-\rho_{k_0}>\rho_{j_0}$. Let $i_1,j_1,k_1$ be defined by
\begin{eqnarray*}
\rho_{i_1} &=& \max \{2\rho_{j_0}-\rho_{k_0}, \rho_{i_0}\} \\
\rho_{j_1} &=& \min \{2\rho_{j_0}-\rho_{k_0}, \rho_{i_0}\} \\
\rho_{k_1} &=& \rho_{j_0}
\end{eqnarray*} 
thus $m=\rho_{i_1}+\rho_{j_1}-\rho_{k_1}$, with $i_1\ge j_1> k_1$ and $i_1\ge i_0,j_1\ge j_0, k_1>k_0$. If $i_1=j_1$, 
then condition b) implies that $m\in S$; otherwise we can repeat the reasoning, obtaining three increasing sequences 
of integers $(i_t), (j_t), (k_t)$, such that
$$
m = \rho_{i_t}+\rho_{j_t}-\rho_{k_t} \\
$$
with $i_t\ge j_t\ge k_t$. There are two possibilities: if there exists an index $h$ such that $i_h=j_h$ or $j_h=k_h$, 
then $m\in S$; otherwise, if $i_t>j_t>k_t$ for all $t$, then, by construction, the sequence $(j_t)$ 
is strictly increasing, so there exists an index $h$ such that $j_h\ge r$, and again we get $m\in S$. 
\end{proof}

\begin{ex}
Let $\scrx$ be the Klein quartic, that is, the curve of homogeneous equation $X^3Y+Y^3Z+Z^3X=0$. 
Let $Q$ be the point at infinity $Q=(1:0:0)$ on $\scrx$. The Weierstrass semigroup of $Q$ is easily seen 
to be $S=\{ 0,3,5,6,7,\cdots\}$. Thus $S$ is an Arf semigroup. 
\end{ex}

\begin{ex} \label{ej}
Let us  consider the tower of function fields $({\scrt}_{n})$ over $\ff_{q^2}$, 
where ${\scrt}_{1}=\ff_{q^{2}}(x_{1})$ and for $n\geq 2$, ${\scrt}_{n}$ is obtained from ${\scrt}_{n-1}$ 
by adjoining a new element $x_{n}$ satisfying the equation
$$ 
x_{n}^{q}+x_{n}=\frac{x_{n-1}^{q}}{x_{n-1}^{q-1}+1}.
$$
This tower was introduced by  Garcia and Stichtenoth in \cite{GS}
(following some previous work of Feng, Rao and Pellikaan)
and it attains the Drinfeld-Vl\u{a}du\c{t} bound. Thus, codes coming from this tower have great interest. 
However, the study of these codes is turning out to be very hard. Some steps in this direction are given 
by H\o holdt and Voss \cite{HV}, Pellikaan, Stichtenoth and Torres \cite{PST}, and Chen \cite{HC}, \cite{HC2}.

Let $Q_{n}$ be the rational place on ${\scrt}_{n}$ that is the unique pole of $x_{1}\,$. 
It is known (see \cite{PST}) that the Weierstrass semigroups $S_{n}$ of ${\scrt}_{n}$ 
at $Q_{n}$ are as follows: $S_1=\no$, and for $n\ge 2$,
$$ 
S_n = q \cdot S_{n-1} \cup \{ m\in\numbers \ | \ m\ge c_n\}
$$
where 
$$
c_{m}=\left\{\begin{array}{ll}
q^{n}-q^{\frac{n+1}{2}} & \mbox{if $n$ is odd;} \\
q^{n}-q^{\frac{n}{2}} & \mbox{if $n$ is even,} \end{array}\right. 
$$
thus, it is easy to see by induction that all of them are Arf semigroups. 

Usually, the codes constructed from this tower are one-point algebraic geometry codes, $C_{\Omega}(\scrp, \rho_lQ_{n})$. 
The minimum distance of some  of these  has been bounded by Chen in \cite{HC},\cite{HC2} where he gives codes on all 
members after level 4 of the family of curves, having true minimum distance greater than the order bound, and uses 
these results to get a sequence of codes giving an improvement on the Tsfasman-Vl\u{a}du\c{t}-Zink bound. 
However note that in those papers, neither the true minimum distance nor the order bound are computed. 
\end{ex}

At the moment, the order bound is already computed for some types of semigroups, including telescopic semigroups 
and semigroups generated by two elements (see \cite{HPvL}). Some results are known for symmetric semigroups 
(see \cite{HPvL} and \cite{Farran}). Remark that Arf semigroups do not lie in these types, because
they are, in general, not symmetric (that is, $c<2g$). The only exception are hyperelliptic semigroups.

\begin{ex}
Let $\scrx$ be an hyperelliptic curve and let $Q$ be a rational hyperelliptic point on $\scrx$. 
The Weierstrass semigroup of $Q$ is hyperelliptic, that is,  $S=\langle 2,t\rangle$, for some odd integer 
$t\ge 3$ (if $t=3$ the semigroup is often called {\em elliptic}). Hyperelliptic semigroups are also Arf semigroups. 
In fact, if $k\le j\le i<r$, then $\rho_i+\rho_j-\rho_k\in 2\numbers \subseteq S$.
\end{ex}

\begin{prop}
The only Arf symmetric semigroups are hyperelliptic semigroups.
\end{prop}
\begin{proof} 
As seen before,
every hyperelliptic semigroup is an Arf semigroup. Conversely, if $S$ is Arf and $\rho\in S, \rho<c$, 
then $\rho+1\not\in S$, because otherwise we have $2(\rho+1)-\rho=\rho+2\in S$, and in the same way 
$\rho+3,\rho+4,\cdots \in S$,  contradicting the fact that $\rho<c$. Thus two consecutive integers in the interval 
$[0,c]$ cannot be  both poles. If $S$ is symmetric, the same happens for gaps (if $l,l+1$ are gaps, then 
$c-l-2,c-l-1$ are poles). Since 0 is always a pole, we get $[0,c]\cap S= [0,c]\cap 2\numbers$ and $S$ is hyperelliptic.
\end{proof}

\section{Computing A-sets in Arf semigroups}

Let  $C_{\Omega}(\scrp,\rho_lQ)$ be an one-point algebraic geometry code. Let $S$ be the Weierstrass semigroup 
at the point $Q$. In the previous sections we have seen how the computation of both the order bound on the 
minimum distance of $C_{\Omega}(\scrp,\rho_lQ)$ and the dimension of the improved codes  related to it, 
involves computations concerning only the semigroup $S$. More precisely, it requires the knowledge of the 
A-sets $A[\rho]$. This study is often difficult for general semigroups.  In this section we shall show that 
the study of the structure and cardinality of the $A[\rho]$'s is rather simple for Arf semigroups.
In order to simplify the exposition, in what follows we shall assume $S\neq\no$. 

For $\rho\in S$, let $j$ be maximum such that $\{ \rho_1,\cdots,\rho_j\} \subseteq A[\rho]$. Then
$$
A[\rho]=\{ \rho_1,\cdots,\rho_j,\rho-\rho_1,\cdots,\rho-\rho_j\}
$$  
because obviously
$\{ \rho_1,\cdots,\rho_j,\rho-\rho_1,\cdots,\rho-\rho_j\} \subseteq A[\rho]$, and conversely, 
if $\rho_k\in A[\rho]$ with $k>j$, we have $\rho-\rho_k=\rho_i$ with $i\le j$, since otherwise 
$\rho-\rho_{j+1}=\rho_i+\rho_k-\rho_{j+1}\in S$, contradicting the choice of $j$.

However, note that the fact $A[\rho]=\{ \rho_1,\cdots,\rho_j,\rho-\rho_1,\cdots,\rho-\rho_j\}$ does not imply  
$\# A[\rho]=2j$, since the set $\{ \rho_1,\cdots,\rho_j,\rho-\rho_1,\cdots,\rho-\rho_j\}$ can contain many 
repeated elements. Thus we define for $\rho\in S$,
\begin{eqnarray*}
\alpha(\rho) &=& \max \{ j \ | \ \rho_1,\cdots,\rho_j\in A[\rho]\} \\
\beta(\rho)  &=& \max \{ j \ | \ \rho_1,\cdots,\rho_j\in A[\rho], \rho_j\le \rho-\rho_j\} \\
 &=& \max \{ j \ | \ \rho_j\in A[\rho], 2\rho_j\le \rho\} .
\end{eqnarray*}
Then $\alpha(\rho)\ge \beta(\rho)$ and we have
$$
A[\rho]=\{ \rho_1,\cdots,\rho_{\beta(\rho)},\rho-\rho_1,\cdots,\rho-\rho_{\beta(\rho)}\}
$$
with
$$
\# A[\rho]=\left\{ \begin{array}{ll}
2\beta(\rho)-1 & \mbox{if $2\rho_{\beta(\rho)}=\rho$;} \\ 
2\beta(\rho)   & \mbox{if $2\rho_{\beta(\rho)}\neq\rho$.} \\
\end{array} \right.
$$ 
In particular, $\# A[\rho]$ is odd if and only if $\rho\in 2S$. The same happens for general semigroups as 
we shall prove later on. Now, let us see what can be said about the numbers $\alpha(\rho)$ and $\beta(\rho)$. 
Let us begin with the case that $\# A[\rho]$ is odd.

\begin{prop}
If $S$ is Arf, then for every $\rho_i\in S$ we have $\beta(2\rho_i)=i$ and consequently $\# A[2\rho_i]=2i-1$.
\end{prop}
\begin{proof}
If $S$ is Arf, then for every $k\le i$ we have $2\rho_i-\rho_k\in S$, so $\{ \rho_1,\cdots,\rho_i\}\subseteq A[2\rho_i]$. 
If $\beta(2\rho_i)>i$ then there exist $j,k>i$ such that $\rho_j+\rho_k=2\rho_i$ what is impossible.
\end{proof}

For poles $\rho\in S\setminus 2S$, we cannot give, in general, an explicit expression for $\beta(\rho)$. 
However we can give some bounds which are enough for our main purposes.

For a positive integer $i$, let $p_i=c+\rho_{i+1}-1$. Clearly $p_i\ge c$ for all $i$, so it is a pole number. 
Furthermore $p_i=\rho_r+\rho_{i+1}-1=\rho_{r+\rho_{i+1}-1}$. In particular, for $i\ge r-1$ we can write $i=(r-1)+t$ 
with $t\ge 0$, and thenwe have $\rho_{i+1}=\rho_{r+t}=c+t$, hence
$p_i=2c+t-1=\rho_{c+i}$.

\begin{prop} \label{>pi}
Let $S$ be a numerical semigroup (not necessarily Arf). If $\rho\in S$ and $\rho>p_{i-1}$, 
then $\{ \rho_1,\cdots,\rho_i\}\subseteq A[\rho]$. If furthermore $i<r$, then $\# A[\rho]\ge 2i$.
\end{prop}
\begin{proof}
If $\rho>p_{i-1}$, then for $j=1,\cdots,i$, we have $\rho-\rho_j\ge c+\rho_i-\rho_j\ge c$. Thus  $\rho-\rho_j\in S$ 
and $\rho_j\in A[\rho]$, hence $\{ \rho_1,\cdots,\rho_i,\rho-\rho_1,\cdots,\rho-\rho_i\} \subseteq A[\rho]$.
If furthermore $i<r$, then $\rho_i<c$ so $2\rho_i\le p_{i-1}<\rho$ and $\rho_i<\rho-\rho_i$. Thus, 
all the elements in the set $\{\rho_1,\cdots,\rho_i,\rho-\rho_1,\cdots,\rho-\rho_i\}$ are distinct. 
\end{proof}

\begin{rem} \label{>=pi}
If $S$ is Arf, proposition \ref{>pi} means that $\alpha(\rho)\ge i$ provided that $\rho>p_{i-1}$. 
Furthermore, if $i<r$ then also $\beta(\rho)\ge i$.
\end{rem}

\begin{prop} \label{=pi}
If $S$ is an Arf semigroup, then $\alpha(p_i)=i$. Furthermore, if $i<r$ then $\beta(p_i)=i$ and $\# A[p_i]=2i$.
\end{prop}
\begin{proof}
Since $p_i-\rho_{i+1}=c-1\not\in S$, then $\rho_{i+1}\not\in A[p_i]$ and $\alpha(p_i)\le i$. 
The conclusion follows from remark \ref{>=pi} and the fact that $2\rho_i<p_i$ for $i<r$. 
\end{proof}

As said before, when $i\ge r-1$ the sequence $(p_i)$ runs over all poles $\rho\ge p_{r-1}=2c-1$, 
that is, for $j\ge c+r-1$ we have $\rho_j=p_{j-c}$. Since $\alpha(p_{j-c})=j-c$, we obtain
$$
A[\rho_j]=\{\rho_1,\cdots,\rho_{j-c},\rho_j-\rho_1,\cdots,\rho_j-\rho_{j-c}\}.
$$
If $j\ge c+r$, then $\rho_j-\rho_s>\rho_{j-c}$ if and only if $s\le r-1$, and we obtain then expression 
$$
A[\rho_j]=\{\rho_1,\cdots,\rho_{j-c},\rho_j-\rho_1,\cdots,\rho_j-\rho_{r-1}\}
$$
without repeated elements. In particular we get the well known result

\begin{prop} \label{lgrande}
For $j\ge c+r$, we have $\# A[\rho_j]=j-g$.
\end{prop}

As a particular case of this proposition, we have $\# A[\rho_{c+r}]=c+r-g= 2r-1$. 
This number is an upper bound for the cardinalities $\# A[\rho_j]$ when $j\le c+r$.

\begin{prop} \label{cotaalta}
If $j<c+r$, then $\beta(\rho_j)\le r-1$ and $\# A[\rho_j]\le 2r-2<\# A[\rho_{c+r}]$.
\end{prop}
\begin{proof}
It suffices to show that $\beta(\rho_j)\le r-1$. Otherwise, if $\beta(\rho_j) \ge r$ 
for some $j<c+r$, then we have $2\rho_r\le\rho_j\le\rho_{c+r-1}$, what leads to $2c\le 2c-1$.
\end{proof}

\begin{cor}
If $\rho_{r-1}+r\le j<c+r$, then $\beta(\rho_j)=r-1$.
\end{cor}
\begin{proof}
If $\rho_{r-1}+r\le j$ then $p_{r-2}<\rho_j$, and the result follows from remark \ref{>=pi} and proposition \ref{cotaalta}.
\end{proof}

\section{The order bound on the minimum distance}

Keeping the notations as in the introduction, let $C_l=C_{\Omega}(\scrp,\rho_lQ)$ be an algebraic geometry code 
arising from a function field $F$, defined by the system of parity checks $\bh_1,\cdots,\bh_l$, that is 
$$
C_l=\{ \bc\in\fq^n \ | \  \bc\cdot\bh_i=0 \mbox{ for all } i=1,\cdots,l\} .
$$  
The dimension of $C_l$ is lower bounded by $n-l$, with equality when all the checks $\bh_i$ are independent $i=1,\cdots,l$. 
It can be shown that this happens if $\rho_l<n$. The minimum distance of $C_l$ is lower bounded by the {\em Goppa bound},  
$d_{G}(l)=l+1-g$. A better bound on the minimum distance is the
{\em order bound} (or {\em Feng-Rao bound}), given by
$$
d_{ORD}(l)=\min\{\# A[\rho] \ | \ \rho\ge\rho_{l+1}\} .
$$
The order bound is always better than the Goppa bound (in fact, it has been proved to be sharp for a number of codes, 
see \cite{HPvL}, but not always, see \cite{HC}). However, it is usually difficult to compute.When the semigroup 
$S$ is Arf, the results obtained in section 2 provide very quickly the order bound of $C_l$ for all $l$.

\begin{thm} \label{order1} 
Let $S$ be an Arf semigroup of genus $g$ and let $c=\rho_r$ be the conductor of $S$. For $i=1,\cdots,r-1$, 
let $l_i=r+\rho_{i+1}-2$. In addition, let $l_0=0$. Then, for any positive integer $l$, we have:
\newline a) if $l_{i-1}<l\le l_i\le l_{r-1}$, then $d_{ORD}(l)=2i$;
\newline b) if $c+r-2=l_{r-1}\le l$, then $d_{ORD}(l)=d_G(l)=l+1-g$.
\end{thm}
\begin{proof}
Since $l_{r-1}=r+c-2$, part b) follows from proposition \ref{lgrande}. To prove a), let us first note that 
for $i=1,\cdots,r-1$, it holds that $p_i=\rho_{r+\rho_{i+1}-1}=\rho_{l_i+1}$. Thus, if $l_{i-1}<l\le l_i$, 
we have $p_{i-1}<\rho_{l+1}\le p_i$, hence, according to propositions \ref{>pi}, \ref{=pi} and \ref{cotaalta}, we have
$$
d_{ORD}(l)=\min\{\# A[\rho] \ | \ \rho\ge\rho_{l+1}\}=\# A[p_i]=2i
$$
and the proof is complete. 
\end{proof}

For some particular types of Arf semigroups we can still give more explicit formulas. For example, 
while studying the redundancy of improved codes coming from the tower in example \ref{ej}, 
Pellikaan and Torres introduce in \cite{PT} the following: 

\begin{defi}\label{defInd}
A sequence $(H_n)$ of semigroups is called {\em inductive} if there exist sequences $(a_n)$ and $(b_n)$ 
of positive integers such that $H_1=\no$ and for $n>1$, $H_n=a_nH_{n-1}\cup\{ m\in \no \ | \ m\ge a_nb_{n-1}\}$. 
A semigroup is called {\em inductive} if it is a member of an inductive sequence.
\end{defi}

The Weierstrass semigroups obtained from the tower of function fields of example \ref{ej} at the points $Q_n$ 
are obviously inductive. Notice that $H_{n}=H_{n-1}$ if $a_{n}=1$. Thus, we can assume that $a_{n}\geq 2$ for $n\geq 2$, 
and hence the sequence $b_{n}$ is super-increasing. For $n\geq 2$ the conductor of $H_{n}$ is obviously 
$c_{n}=a_{n}b_{n-1}$\/. Since $\no$ is inductive, with the aid of the following result one easily proves 
by induction that any inductive semigroup is Arf. 

\begin{lemma}\label{lemInd}
Let $S$ be an Arf semigroup and take arbitrary positive integers $a,R$. 
Then $\overline{S}=aS \cup\{ m\in\no \ | \ m\ge R\}$ is an Arf semigroup. 
\end{lemma}
\begin{proof}
Let $\overline{\rho_i},\overline{\rho_j},\overline{\rho_k}\in\overline{S}$, $i\ge j\ge k$, be three poles 
smaller than the conductor of $\overline{S}$ (hence $R>\overline{\rho_i}\ge\overline{\rho_j}\ge\overline{\rho_k}$). 
There exist poles $\rho_{\alpha},\rho_{\beta},\rho_{\gamma}\in S$ such that $\alpha\ge\beta\ge\gamma$ and 
$\overline{\rho_i}=a\rho_{\alpha},
\overline{\rho_j}=a\rho_{\beta},
\overline{\rho_k}=a\rho_{\gamma}$. Then, since 
$\overline{\rho_i}+\overline{\rho_j}-\overline{\rho_k}=a(
\rho_{\alpha}+\rho_{\beta}-\rho_{\gamma})\in aS\subseteq\overline{S}$,
the result follows from the fact that $S$ is Arf. 
\end{proof}

As a consequence, given an inductive sequence of semigroups, $(H_n)$, in order to determine the order bound for $H_n$ 
one can describe inductively the intervals where such bound changes, according to the results of the previous section.  
Let $c^{(n)}, r^{(n)}, \rho_i^{(n)} (i=1,2,\cdots)$ and $l_i^{(n)} (i=0,\cdots,r^{(n)}-1)$, be the corresponding elements 
and parameters of $H_n$. In addition, let $g^{(n)}$ be the genus of $H_n$ and denote 
$\lambda^{(n)}=b_n-c^{(n)}, \lambda^{(0)}=1, L^{(n)}=\lambda^{(0)}+\cdots+\lambda^{(n)}$.

\begin{prop} 
With the above notations, the following holds: \newline
a) $c^{(n)}=a_nb_{n-1}$; \newline
b) $r^{(n)}=L^{(n-1)}$; \newline
c) for $i=1,\cdots,r^{(n)}$, we have $\rho_i^{(n)}= a_n\rho_i^{(n-1)}$, and hence \newline
c.1) for $i=1,\cdots,r^{(n-1)}-1$ one has $\rho_{i+1}^{(n)}=a_{n}\rho_{i+1}^{(n-1)}$, and thus  
$l_{i}^{(n)}=l_{i}^{(n-1)}+\lambda^{(n-1)}+(a_{n}-1)\rho_{i+1}^{(n-1)}$;
\newline
c.2) for $i=r^{(n-1)}+1,\cdots,r^{(n)}$ one has $\rho_{i}^{(n)}=a_{n}(c^{(n-1)}+i-r^{(n-1)} )$, and thus
$l_{i}^{(n)}=r^{(n-1)}+\lambda^{(n-1)}-2+a_{n}(c^{(n-1)}+i+1-r^{(n-1)})= r^{(n)}-2+a_{n}(c^{(n-1)}+i+1-r^{(n-1)} )$; 
\newline
d) $g^{(n)}=a_nb_{n-1}-L^{(n-1)}+1$.
\end{prop}

The proof of this result is left to the reader.

There is a nice alternative description of the semigroup $H_{n}$ which allows us 
to compute in another way the intervals where the order bound is constant. 
Namely, such intervals are described in an iterative way, instead of recursively. 
In fact, for $k=1,\cdots,n-1$, denote $A_{k}^{(n)}=\prod_{i=k+1}^{n}a_{i}$.  
Then $H_{n}$ can be described as follows:
$\rho_{1}^{(n)}=0$; the following $\lambda^{(1)}$ poles are obtained by summing $A_{1}^{(n)}$ 
to the previous one; the following $\lambda^{(2)}$ poles are obtained by summing $A_{2}^{(n)}$ 
to the previous one; and so on until we reach $c^{(n)}$, and then we sum $1$ each time. This description of $H_n$ 
will be called $[\star]$. It allows us to list the poles $\rho_{i}^{(n)}$ and the numbers $l_{i}^{(n)}$ for $H_{n}$ 
by means of the following  

\begin{prop} \label{torre}
With the above notations, if $L^{(k)}<i\le L^{(k+1)}$ and $\lambda^{(k+1)}>0$ then 
$$\rho_{i}^{(n)}=\rho_{L^{(k)}}^{(n)}+(i-L^{(k)})A_{k+1}^{(n)}$$
and hence 
$$
l_{i-1}^{(n)}=L^{(n-1)}-2+\rho_{i}^{(n)}.
$$
\end{prop}

\begin{ex}
Consider again the tower of function fields $({\scrt}_{n})$ given in example \ref{ej}. 
Since the semigroups $S_{n}$ are inductive, one can apply the above results 
to compute the order bound. In this way, we obtain: 

\begin{itemize}
\item $a_{n}=q$ for $n\geq 2$, 
\item $A_{k}^{(n)}=q^{n-k}$ for $1\leq k\leq n-1$, 
\item $\lambda^{(2i-1)}=q^{i-1}(q-1)$ and $\lambda^{(2i)}=0$ for $i\geq 1$, 
\item $L^{(2i-1)}=L^{(2i)}=q^{i}$ for $i\geq 1$, hence $L^{(n)}= q^{\lfloor (n+1)/2 \rfloor}$. 
\end{itemize}
\noindent
By using the description $[\star]$ for $S_{n}$, one easily obtains
$$ 
\rho_{q^{k}}^{(n)}=q^{n-k}(q^{k}-1)
$$
for $k=0,\cdots,\lfloor n/2\rfloor$. Then, if $q^k<i+1\le q^{k+1}$ for some $k$, 
with $0\le k\le \lfloor n/2\rfloor$, from proposition \ref{torre} we get
\begin{eqnarray*}
l_{i}^{(n)}&=& q^{\lfloor\frac{n}{2}\rfloor}-2+(i+1-q^{k})q^{n-k-1}+q^{n-k}(q^{k}-1) \\
&=& q^{\lfloor\frac{n}{2}\rfloor}-2+q^{n-k-1}(q^{k+1}-q^{k}-q+i+1).
\end{eqnarray*}

Since $r^{(n)}=L^{(n-1)}=q^{\lfloor n/2\rfloor}$, this formula provides all the values 
$l_1^{(n)},\cdots, l_{r^{(n)}-1}^{(n)}$, and hence, according to theorem \ref{order1}, 
the order bound for all codes $C_l$ coming from $\scrt_n$.
\end{ex}

\section{The redundancy of improved codes}

By using the same notation as in the previous section, let us consider the algebraic geometry code $C_l$ 
defined by means of the set of checks $\bh_1,\cdots,\bh_l$. For a positive integer $d$ let us consider the set
$$
R_d=\{ i \ | \ \# A[\rho_i]<d\} \sim \{ \rho\in S \ | \ \# A[\rho]<d\}
$$
and the {\em improved geometric Goppa code} $\tilde{C}(d)$ defined as
$$
\tilde{C}(d)=\{ \bc\in\fq^n \ | \ \bc\cdot\bh_i=0 \mbox{ for all } i\in R_d\}
$$
(see \cite{FR,HPvL,PT}). The minimum distance of $\tilde{C}(d)$ is at least $d$. Furthermore, if $d=d_{ORD}(l)$, 
then $C_l\subseteq\tilde{C} (d)$ (this is the reason of the term \lq improved'). Thus, a natural question is 
to compute the improvement on the dimension, $\dim \tilde{C}(d)-\dim C_l$. It is well known (see \cite{HPvL}) 
that $\dim C_l\ge n-l$, with equality if $\rho_l<n$; on the other hand, from its definition, it follows that 
$\dim \tilde{C}(d)\ge n-\# R_d$. It is easy to see that when $2c\le n$, then  for $d$ in the range
$1\le d\le 2r-2$ (where one can hope an improvement on the dimension) all
the checks $\bh_i$ are independent, and thus we have equality, $\dim \tilde{C}(d)=n-\# R_d$. In this section we shall 
compute the sequence $(\# R_d)$ when the semigroup $S$ is Arf. This result, together with theorem \ref{order1}, 
allows us the computation of the improvement on the dimension.

The sequence $(\# R_d)$ has been already treated in the paper \cite{PT} by Pellikaan and Torres.
They show that for every semigroup $S$, one has
\begin{eqnarray*}
\# R_d &=& d+g-1 \; \mbox{ if $d\ge 2r-1$;} \\
\# R_{2r-2} &=& \rho_{r-1}+r-1 
\end{eqnarray*}
what can be written as
$$
\# R_d=\rho_{\lceil\frac d2 \rceil}+\lfloor\frac d2 \rfloor
$$
provided that $d\ge 2r-2$. In order to simplify the exposition, the semigroups verifying the above formula for all 
$d\ge 1$ will be called {\em stable}. In the paper \cite{PT}, the authors propose the characterization of stable 
semigroups as an open problem, and they show that inductive semigroups are stable.

In this section, we shall prove that stable semigroups are precisely Arf semigroups. To prove this, we introduce 
some notation. For $d$ a positive integer, let $S_d=\{ \rho\in S \ | \ \# A[\rho]=d\}$ (so $\# R_{d+1}=\# R_d+\# S_d$). 
In order to characterize stable semigroups it is enough to consider values of $d$ in the range that $1\le d \le 2r-3$, 
and then, a semigroup $S$ is stable if and only if for every odd integer $d$, $1\le d\le 2r-3$, if we write $d=2t+1$, 
then it holds that
\begin{eqnarray*}
\# S_d &=& 1 \\
\# R_d &=& \rho_{t+1}+t.
\end{eqnarray*}

\begin{lemma} \label{2i-1}
Let $S$ be a semigroup and let $\rho\in S$. Then $\# A[\rho]$ is odd if and only if $\rho\in 2S$. 
In this case, if $\rho=2\rho_i$ then $\# A[\rho]\le 2i-1$.
\end{lemma} 
\begin{proof}
For every $p\in A[\rho]$ we have $p'=\rho-p\in A[\rho]$, hence $\# A[\rho]$ is even unless there exists a 
(unique a fortiori) pole $p\in A[\rho]$ such that $p=p'=\rho-p$, that is, $\rho\in 2S$. In this case, 
if $\rho=2\rho_i$, then for every $p,p'\in A[\rho]$ with $p+p'=\rho=2\rho_i$, then either $p\le\rho_i$ 
or $p'\le \rho_i$, and hence $\# A[\rho]\le 2i-1$. 
\end{proof}

\begin{prop} \label{casiArf}
Let $S$ be a semigroup. The following statements are equivalent:
\newline a) $\# A[2\rho_i]=2i-1$ for all $\rho_i\in S$;
\newline b) $\# S_d=1$    for all $d$ odd;
\newline c) $S$ is Arf.
\end{prop}
\begin{proof}
According to lemma \ref{2i-1}, we have $\# S_d=1$ for all $d$ odd if and only if $\# A[2\rho_i]=2i-1$ for all  
$\rho_i$, and this happens if and only if 
$$
A[2\rho_i]=\{ \rho_1,\cdots,\rho_i,2\rho_i-\rho_1,\cdots,2\rho_i-\rho_i \}
$$
that is, if and only if $2\rho_i-\rho_j\in S$ for all $i, j$ with $i\ge j$. This is equivalent to $S$ 
being Arf according to proposition \ref{debil}.
\end{proof}

Thus, all stable semigroups are Arf.
In the sequel we shall prove that Arf semigroups are stable. To that end it suffices to show that for all $d$ 
odd in the range $1\le d \le 2r-3$, if $d=2t+1$, then  $\# R_d=\rho_{t+1}+t$.

\begin{lemma}
Let $S$ be an Arf semigroup and let $d$ be as above. Then $R_d\subseteq [0,p_t]\cap S$.
\end{lemma}
\begin{proof}
If $\# A[\rho]<d\le 2r-3$, then proposition \ref{cotaalta} implies $\rho\le p_{r-1}$. 
Thus the result follows from proposition \ref{>pi}. 
\end{proof}

For $\rho\in [0,p_t]\cap S$, it holds that $\rho\in R_d$ if and only if $\beta(\rho)\le t$, so we get the following

\begin{lemma} \label{44}
Let $S$ and $d$ be as in the previous lemma. Then
$$
\{\rho\in [0,p_t]\cap S \ | \ \beta(\rho)\ge t+1 \} =
\{\rho_{t+1}+\rho_{t+1},\cdots,\rho_{t+1}+\rho_{r-1}\}. 
$$
\end{lemma}
\begin{proof}
If $\rho\in [0,p_t]\cap S$ is such that $\alpha (\rho)\ge t+1$, then we have $\rho=\rho_{t+1}+\rho_i$ for some $i$. 
Since $2\rho_{t+1}\le \rho$, it holds that $i\ge t+1$. On the other hand, since $\rho\in [0, p_t]$ and $\beta(p_t)=t$, 
we have $\rho< p_t$ and $i\le r-1$. Hence
$\{\rho\in [0,p_t]\cap S \ | \ \beta(\rho)\ge t+1 \} \subseteq\{\rho_{t+1}+\rho_{t+1},\cdots,\rho_{t+1}+\rho_{r-1}\}$. 
The converse is clear. 
\end{proof}

Finally we have the following

\begin{thm}
Let $S$ be a semigroup. The following statements are equivalent:
\newline a) $S$ is Arf;
\newline b) for every positive integer $d$, we have
$\# R_d=\rho_{\lceil\frac d2 \rceil}+\lfloor\frac d2 \rfloor$.
\end{thm} 
\begin{proof}
If b) holds then $\# S_d=1$ for all $d$ odd and $S$ is an Arf semigroup as we have seen in proposition \ref{casiArf}. 
Conversely, assume $S$ is Arf and let $d$ be an odd integer with $1\le d\le 2r-3$. According to proposition \ref{casiArf} 
we have $\# S_d=1$. Now if we write $d=2t+1$, then,
according to lemma \ref{44}, we have
$$
\# R_d=\# ([0,p_t]\cap S)-\# \{\rho_{t+1}+\rho_{t+1},\cdots,\rho_{t+1}+\rho_{r-1}\}.
$$
Since $p_t=\rho_r+\rho_{t+1}-1=\rho_{r+\rho_{t+1}-1}$, we obtain $\# ([0,p_t]\cap S)=r+\rho_{t+1}-1$.
Thus, $\# R_d=(r+\rho_{t+1}-1)-(r-t-1)= \rho_{t+1}+t$ and $S$ verifies b). 
\end{proof}

We are now able to compare the dimension of  the codes $C_l$ and $\tilde{C}(d)$.

\begin{prop}
Let $S$ be an Arf semigroup. For a positive integer $l$ let us consider
the codes $C_l$ and $\tilde{C}(d)$, where $d=d_{ORD}(l)$.
Let $l_0,\cdots,l_{r-1}$ be as in theorem \ref{order1}. 
\newline a) If $l_{i-1}<l\le l_i$, with $i\le r-1$ and $2c\le n$, then $\dim\tilde{C}(d)-\dim C_l=l-\rho_i-i$.
\newline b) If $l>l_{r-1}=c+r-2$, then  $C_l=\tilde{C}(d)$.
\end{prop}
\begin{proof}
If $2c\le n$, then (see \cite{PT}) all the checks $\bh_i$ in $C_l$ and $\tilde{C}(d)$ are independent, 
so $\dim C_l=n-l$, $\dim \tilde{C}(d)=n-\# R_d$ and $\dim \tilde{C}(d)-\dim C_l=l-\# R_d$. Now, if $l_{i-1}<l\le l_i$, 
then $d=2i$ and $\# R_d=\rho_i+i$. Thus $l-\# R_d=l-\rho_i-i$. This proves a). If $l\ge c+r-1$, then, according 
to theorem \ref{order1}, we have $d=l+1-g\ge 2r-1$. Thus $\# R_d=d+1-g=l$,  hence $R_d=\{\rho_1,\cdots,\rho_l \}$ and 
$C_l=\tilde{C}(d)$.
\end{proof}

\vspace{1.cm}
\Hline
\vspace{-1.cm}


\begin{thebibliography}{99}

\bibitem{Ar} C. Arf, \lq\lq Une interpretation alg\'{e}brique de la suite des ordres de 
multiplicit\'{e} \ d'une branche alg\'{e}brique", {\em Proc. London Math. Soc.} vol. 50, pp. 256-287, 1949. 

\bibitem{Farran} A. Campillo and J. I. Farr\'{a}n, \lq\lq Computing Weierstrass semigroups and 
the Feng-Rao distance from singular plane models", to appear in {\em Finite fields and their applications}.  

\bibitem{HC} H. Chen, \lq\lq Codes on Garcia-Stichtenoth curves with true distance greater than  Feng-Rao distance", 
{\em IEEE Trans. Inform. Theory}, vol IT-45, pp. 706-709, March 99. 

\bibitem{HC2} H. Chen, \lq\lq On the number of correctable errors of the Feng-Rao decoding algorithm for AG codes", 
{\em IEEE Trans. Inform. Theory}, vol IT-45, pp. 1709-1712, July 99. 

\bibitem{FRB} G. L. Feng and T. R. N. Rao, \lq\lq A simple approach for construction of algebraic geometry codes 
from affine plane curves", {\em IEEE Trans. Inform. Theory}, vol. IT-40, pp. 1003-1012, July 1994. 

\bibitem{FR} G. L. Feng and T. R. N. Rao, \lq\lq Improved Geometric Goppa codes, Part I: Basic Theory", 
{\em IEEE Trans. Inform. Theory}, vol. IT-41, pp. 1678-1693, Nov. 1995. 

\bibitem{GS} A. Garc\'{\i}a and H. Stichtenoth, \lq\lq On the asymptotic behaviour of some towers of function fields 
over finite fields", {\em J. Number Theory}, vol. 61, pp. 248-273, 1996. 

\bibitem{HV} T. H\o holdt and C. Voss, \lq\lq An explicit construction of a sequence of codes attaining 
the Tfasman-Vl\u{a}du\c{t}-Zink bound: the first steps", {\em IEEE Trans. Inform. Theory}, 
vol. IT-43, pp. 128-135, Jan. 1997. 

\bibitem{HPvL} T. H\o holdt, J. H. van Lint and R. Pellikaan, \lq\lq Algebraic Geometry codes", 
in {\em Handbook of Coding Theory}, V. Pless, W. C. Huffman and R. A. Brualdi, Eds, pp. 871-961 (vol 1), 
Elsevier, Amsterdam, 1998. 

\bibitem{Li} J. Lipman, \lq\lq Stable ideal and Arf rings", {\em Amer. J. Math.} vol. 97, pp. 791-813, 1975. 

\bibitem{PST} R. Pellikaan, H. Stichtenoth and F. Torres, \lq\lq Weierstrass semigroups in an asymptotically good 
tower of function fields", to appear in {\em Finite fields and their applications}. 

\bibitem{PT} R. Pellikaan and F. Torres, \lq\lq On Weierstrass semigroups and the redundancy 
of improved geometric Goppa codes", preprint, 1998. 

\end{thebibliography}
\end{document}